\documentclass[lettersize,journal]{IEEEtran}
\usepackage{amsmath,amssymb,amsfonts,amsthm}
\usepackage{algpseudocode}
\usepackage{algorithm}
\usepackage{array}
\usepackage[caption=false,font=normalsize,labelfont=sf,textfont=sf]{subfig}
\usepackage{textcomp}
\usepackage{stfloats}
\usepackage{url}
\usepackage{verbatim}
\usepackage{graphicx}
\usepackage{cite}
\usepackage{bm}
\usepackage{booktabs}
\usepackage{dsfont}
\usepackage{textcomp}
\usepackage{enumitem}
\usepackage{multirow}
\usepackage{hyperref}

\usepackage[markup=default,commandnameprefix=always,authormarkup=none]{changes}

\newtheorem{theorem}{Theorem}
\newtheorem{lemma}{Lemma}

\newtheorem{definition}{Definition}

\newtheorem{remark}{Remark}

\begin{document}

\title{False Data-Injection Attack Detection in Cyber-Physical Systems: A Wasserstein Distributionally Robust Reachability Optimization Approach

\author{Yulin Feng, Dapeng Lan, \IEEEmembership{Member, IEEE}, and Chao Shang, \IEEEmembership{Member, IEEE}} 
\thanks{This work was supported by the National Natural Science Foundation of China (Grant No. 62373211), and the Open Research Project of the State Key Laboratory of Industrial Control Technology, China (Grant No. ICT2025B10) \textit{(Corresponding author: Chao Shang)}.}
\thanks{Yulin Feng and Chao Shang are with the Department of Automation, Beijing National Research Center for Information Science and Technology, Tsinghua University, Beijing 100084, China (e-mail: fyl23@mails.tsinghua.edu.cn; c-shang@tsinghua.edu.cn;).}
\thanks{
Dapeng Lan is with the Shenyang Institute of Automation, Chinese Academy of Sciences, Shenyang 110003, China (e-mail: landapeng@sia.cn).}
}
 
\markboth{\tiny This work has been submitted to the IEEE for possible publication. Copyright may be transferred without notice, after which this version may no longer be accessible.}
{Shell \MakeLowercase{\textit{et al.}}: A Sample Article Using IEEEtran.cls for IEEE Journals}


\maketitle

\begin{abstract}
Cyber-physical system (CPS) is the foundational backbone of modern critical infrastructures, so ensuring its security and resilience against cyber-attacks is of pivotal importance. This paper addresses the challenge of designing anomaly detectors for CPS under false-data injection (FDI) attacks and stochastic disturbances governed by unknown probability distribution. By using the Wasserstein ambiguity set, a prevalent data-driven tool in distributionally robust optimization (DRO), we first propose a new security metric to deal with the absence of disturbance distribution. This metric is designed by asymptotic reachability analysis of state deviations caused by stealthy FDI attacks and disturbance in a distributionally robust confidence set. We then formulate the detector design as a DRO problem that optimizes this security metric while controlling the false alarm rate robustly under a set of distributions. This yields a trade-off between robustness to disturbance and performance degradation under stealthy attacks. The resulting design problem turns out to be a challenging semi-infinite program due to the existence of distributionally robust chance constraints. We derive its exact albeit non-convex reformulation and develop an effective solution algorithm based on sequential minimization. Finally, a case study on a simulated three-tank is illustrated to demonstrate the efficiency of our design in robustifying against unknown disturbance distribution.
\end{abstract}

\section{Introduction}
Cyber-Physical Systems (CPSs) are a critical component of modern technological advancements, integrating cyber communication and computation with physical plants. CPSs can be used in various fields, such as autonomous vehicles \cite{santoso2022data}, industrial control systems \cite{zhang2022advancements} and power grid \cite{edib2023cyber}. However, CPSs are often exposed to threats from external cyber-attackers, who may compromise communication networks, manipulate sensor data, or interfere with control signals, leading to severe consequences. Therefore, the security of CPS has gained increasing research attention over the recent decade.

As an essential problem of CPS security, anomaly detection aims at identifying unusual or malicious actions, thereby addressing the underlying challenge of mitigating malicious activities before significant harms are caused \cite{zhang2022design}. The anomaly detection problem is closely related to fault detection of technical processes, which has been a hot spot in the control community \cite{ding2008model,ding2014data}. However, cyber-attacks are more difficult to deal with than faults since attackers can meticulously design the attack mechanism to deceive the detector by exploiting available system knowledge \cite{zhao2020data,zhang2022design,lu2022false,wu2022covert}. Among various types of cyber-attacks, the stealthy false data-injection (FDI) attacks are a critical threat, as they can cause severe performance degradation while bypassing attack alarm\cite{wu2022covert,renganathan2020distributionally,zhang2022design}. Therefore, many tailored detector design schemes have been developed specifically to counter such attacks. In \cite{zhao2020data}, a subspace-based detector based on coding theory was proposed to tackle undetectable sensor FDI attacks. \cite{yang2019multiple} borrowed ideas from the random finite set theory to detect multiple attacks on different sensors. Deep reinforcement learning has also been found useful for detecting FDI attacks under disturbance by formulating the design problem as a partially observable Markov decision process in \cite{liu2022false}. Aiming to accelerate the response speed, \cite{gao2019fusion} designed an optimal weighting fusion criterion to calibrate the threshold under the limited bandwidth. A generalized likelihood ratio-based scheduler was presented in \cite{ding2017likelihood}, which selectively transmits the most informative sensor data to detect potential cyber-attacks under limited communications. 

Uncertain disturbance is widespread in real-world CPS. In previous works on cyber-attack detection, it is frequently assumed that the disturbance is either bounded or governed by a Gaussian distribution to ease analysis and design.
However, once disturbance distribution deviates from Gaussianity, such as exhibiting heavy-tailed or non-stationary characteristics, the accuracy and reliability of cyber-attack detection are inevitably compromised. As an effective technique for managing uncertainty in probability distributions, distributionally robust optimization (DRO) has attracted extensive attention in the field of anomaly detection recently. Instead of assuming the true probability distribution to be precisely known, DRO optimizes the worst-case performance within the ambiguity set composed of all possible distributions, thereby offering a more robust solution than generic methods.  Although various distributionally robust anomaly detectors have already been developed, see e.g. \cite{xue2022integrated,cheng2025improved,shang2021distributionally,shang2022generalized,feng2024distributionally}
While various distributionally robust anomaly detectors have already been developed, see e.g. \cite{xue2022integrated,cheng2025improved,shang2021distributionally,shang2022generalized,feng2024distributionally}, they are primarily oriented towards fault detection tasks. These methods merely aim at maximize detectability, but overlook the performance degradation caused by attacks, which also serves as a essential metric \cite{mo2015performance,li2019performance,liu2024performance}.

This work aims at addressing the FDI attack detection problem in CPS subject to stochastic disturbance whose distribution is unknown. The main contributions of this work can be summarized as follows: 
\begin{itemize}
\item 
We propose a new distributionally robust security metric to evaluate the performance degradation due to stealthy FDI attacks and stochastic disturbance governed by a unknown distribution. This is achieved by using the Wasserstein ambiguity set for data-driven uncertainty description and asymptotic reachability analysis.

\item We formulate FDI attack detector design as a DRO problem that minimizes our proposed performance degradation metric while simultaneously controlling the false alarm rate (FAR), thereby achieving a desirable trade-off between security against attacks and robustness to disturbances.
\item We establish an exact reformulation of the proposed distributionally robust detector design problem, successfully transforming a semi-infinite program into a finite-dimensional problem but involving bilinear matrix inequalities (BMIs). To solve this resulting non-convex problem, we develop a customized and efficient algorithm based on sequential minimization.
\end{itemize}

The remainder of this article is organized as follows.
The basics of CPS and the distributionally robust anomaly detector design are revisited in Section \ref{preli}. Section \ref{mainr} presents our attack detector design scheme, and case study results are given in Section \ref{case}. Section \ref{concl} concludes this article. 

 \textit{Notations:} We use $\mathbb{N}_{i:j}$ to denote the set of integer indexes $\{i,\cdots,j\}$. $\mathcal{N}(\mu,\Sigma)$ denotes a Gaussian distribution with mean $\mu$ and covariance $\Sigma$. $\chi_{m}^2$ denotes a Chi-square distribution with $n$ degrees of freedom and its upper $\varepsilon$-quantile is denoted by $\chi_{m}^2(\varepsilon)$. $I$ denotes an identity matrix of appropriate dimension. For a vector $x$, the weighted $l_2$-norm is denoted by $\lVert x \rVert_{W} =(x^\top W x)^{\frac{1}{2}}$, $\lVert f \rVert$ denotes the Euclidean norm, and  $\lVert f \rVert_{\infty}$ denotes the $l_{\infty}$-norm. For a matrix $X$, $X_{[i:j]}$ denotes the submatrix of $X$ that goes from the $i$th to the $j$th column. $X^\dagger$ denotes its Moore-Penrose inverse, ${\rm tr}\{X\}$ denotes its trace, and $\lVert X\rVert_{2} $ denotes its spectral norm. ${\rm diag}\{X_1,\cdots,X_n\}$ is the block diagonal matrix with diagonal block matrices $X_1,\cdots,X_n$. For a symmetric $X$, $X \succeq(\succ) 0$ indicates that $X$ is positive (semi-)definite. For a discrete-time signal $x(k)$, the concatenated vector is denoted by $x_s(k) = \left[x(k-s+1)^\top \cdots x(k)^\top \right]^\top$.

\begin{figure}
    \centering
    \includegraphics[width=0.95\linewidth]{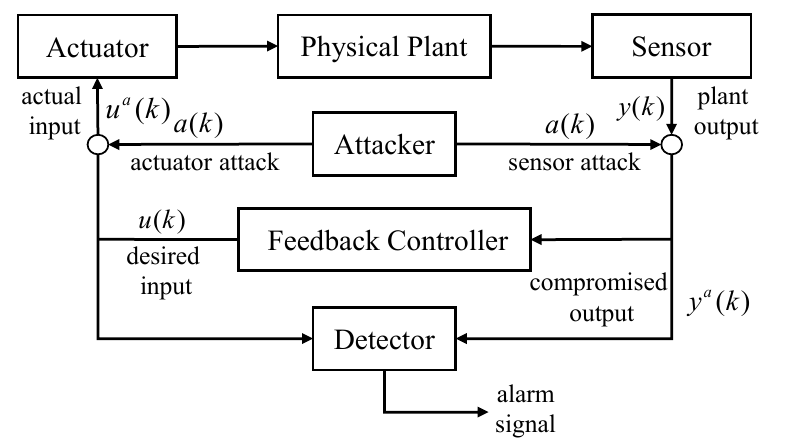}
    \caption{CPS under FDI attacks.}
    \label{fig_CPS}
\end{figure} 

\section{Preliminaries} \label{preli}
\subsection{CPS Configuration}
As described in Fig. \ref{fig_CPS}, a CPS under study consists of a physical plant, a feedback controller, and an attack detector. The physical plant can be described as a stochastic discrete-time linear time-invariant (LTI) system: 
\begin{equation}\begin{cases} {x}(k+1)=A{x}(k)+B {u}^a(k)+B_d d(k)\\
{y}(k)=C {x}(k)+D_dd(k)\end{cases} \label{eq_SS}\end{equation}
where ${x}\in\mathbb{R}^{n_{x}}$, ${y}\in\mathbb{R}^{n_y}$, ${u}^a\in\mathbb{R}^{n_u}$ and $d\in\mathbb{R}^{n_d}$ denote system state, output signal, actual control signal, stochastic disturbance, and additive disturbance, respectively. $A$, $B$, $B_d$, $C$ and $D_d$ are state-space matrices of appropriate dimensions. It is assumed that $(A, C)$ is observable and $(A, B)$ is controllable. 
The sensor and actuator attacks can be described as:
\begin{equation}\begin{cases} 
 {u}^a(k)= u(k)+ B_a a(k)\\
 {y}^a(k)= y(k)+ D_a a(k)
\end{cases} \label{eq_SS_atk}
\end{equation}
where the attack $a(k)\in\mathbb{R}^{n_u+n_y}$ signal injects false data into $u(k)$ and $y(k)$ through channels $B_a$ and $D_a$,  and \( y^a(k) \) represents the compromised output. \( u(k) \) is the desired input produced by the feedback controller:
\begin{equation}\begin{cases} {x}_c(k+1)=A_c{x}_c(k)+B_c[y_{\rm ref}(k)-y^a(k)]\\
{u}(k)=C_c {x}_c(k)+D_c[y_{\rm ref}(k)-y^a(k)]\end{cases} \label{eq_SS_ctrl}\end{equation}
where ${x}_c\in\mathbb{R}^{n_{c}}$ and $y_{\rm ref}(k)\in\mathbb{R}^{n_y}$ are the state of the controller and the output reference. By defining the augmented state ${\bar x}(k)=\begin{bmatrix}
    { x}(k)^\top & { x}_c(k)^\top
\end{bmatrix}^\top$ and combining \eqref{eq_SS}, \eqref{eq_SS_atk}
 and \eqref{eq_SS_ctrl}, one attains the dynamics of the closed-loop system under attack:
\begin{equation} 
\begin{cases}
    \bar {x}(k+1)=\bar A{\bar x}(k)+ \bar B_{\rm r} {y_{\rm ref}(k)}+\bar B_a a(k)+ \bar B_d d(k)\\
    y^a(k) = \bar C{\bar x}(k)+ D_a a(k)+D_d d(k)
\end{cases}\label{eq_SS2}
\end{equation}
where 
\begin{align*}
    &\bar{A}=\begin{bmatrix}
        A-BD_cC&BC_c\\
       -B_cC &A_c
    \end{bmatrix},~\bar B_{\rm r}=\begin{bmatrix}
        BD_c\\B_c
    \end{bmatrix},~\bar C= \begin{bmatrix}
       C\\0
    \end{bmatrix},\\
   & \bar B_a = \begin{bmatrix}
     B(B_ a- D_c D_a)  \\
       -B_c D_a
   \end{bmatrix},~\bar B_{d}=\begin{bmatrix}
        B_d-BD_cD_d\\-B_cD_d
    \end{bmatrix}.
\end{align*}

To detect possible attacks, a parity-space-based residual generator can be constructed. Given order $s \geq n_x$, the parity relation of the desired input $u_s(k)$ and the compromised output $y^a_s(k)$ is expressed as
\begin{equation} \label{eq_parity}
y^a_s(k) = \Gamma_s x(k-s)+ H_{u,s}u_s(k) + H_{d,s}d_s(k) + H_{a,s}a_s(k),
\end{equation}
where $\Gamma_s = \begin{bmatrix}C^\top & (CA)^\top & \dots & \left(CA^{s-1}\right)^\top\end{bmatrix}^\top$ is the extended observability matrix, and
\begin{equation*}
    H_{u,s}=\begin{bmatrix}D&0&0&\cdots&0\\
    CB&D&0&\cdots&0\\
    CAB&CB&D&\cdots&0\\
    \vdots&\ddots&\ddots&\ddots&\vdots\\CA^{s-1}B&\cdots&CAB&CB&D\end{bmatrix},\\
\end{equation*}
is the Toeplitz matrix. Then $H_{d,s}$ and $H_{a,s}$ can be constructed in a similar form of  $H_{u,s}$ with  $\{B,D\}$ replaced by  $\{B_d,D_d\}$ and  $\{BB_a,D_a\}$, respectively. 
As a result, the residual signal $r(k)\in \mathbb{R}^{n_r}$ can be generated as follows \cite{chow1984analytical}:
\begin{equation} 
\begin{aligned}
    r(k)& = P \cdot \Gamma^{\perp}_s \left[  y^a_s(k) - H_{u,s}u_s(k) \right]\\
    &= P[W_dd_s(k)+W_aa_s(k)]
\end{aligned}\label{eq_fd_parity}
\end{equation}
where the disturbance and attack projection matrices $W_d=\Gamma^{\perp}_s H_{d,s}$, $w_a=  \Gamma^{\perp}_s H_{a,s}$, $\Gamma^{\perp}_s\in \mathbb{R}^{(sn_y-n_x)\times sn_y}$ is the orthogonal complement of $\Gamma$ eliminating the impact of the initial state $x(k-s)$ on $r(k)$, and $P$ is the projection matrix with some design freedom. It is obvious that $r(k)$ depends affinely on $d_s(k)$ and $a_s(k)$, and the choice of $P$ shall balance between insensitivity to $d_s(k)$ and sensitivity to $a_s(k)$. It is essential to ensure that the design matrix $P$ has full column rank, thereby not being blind to some not strictly stealthy attack directions $a_s \notin {\rm ker}(W_a)$. Subsequently, one can evaluate the residual through the function $J(r)=\lVert r(k) \rVert^2$ to detect the occurrence of attacks. The alarm logic is expressed as:
\begin{equation*}
    \begin{cases}
    J(r) > J_{{\rm th}} \Rightarrow & \text{Attack alarm} \\
    J(r) \leq J_{{\rm th}} \Rightarrow & \text{No alarm}
    \end{cases}
\end{equation*}
Without loss of generality, the threshold can be set to $ J_{{\rm th}} = 1$, since its design freedom can be accounted for by optimizing the projection matrix $P$. Due to the randomness in $d_s(k)$, unwanted false alarms may be raised under the attack-free condition, and FAR has been widely adopted to evaluate robustness performance of the detector. We use $\xi \in \mathbb{R}^{n_\xi}$ to denote the random variable of the $s$-long augmented disturbance $d_s$, where the time index $k$ is dropped for notational brevity. 
\begin{definition}[FAR \cite{ding2008model}]\label{def_FAR}
Given the threshold $J_{\rm th}=1$ and the evaluation function $J(r)=\lVert r \rVert^2$, the FAR is defined as
\begin{align*}   
{\rm FAR} &= \mathbb{P}_{\xi} \left\{ \lVert r\rVert^2 > 1\mid a_s=0 \right\}\\
&= \mathbb{P} \left\{\lVert P  W_d \xi  \rVert ^2 > 1\right\}.
\end{align*}
\end{definition}
To compute FAR, the probability distribution of $\xi$ has to be known exactly. The Gaussian assumption, i.e., $\xi\sim \mathcal{N}(0,\Sigma_0) $, is widely adopted in present literature. As an example, in the GLRT approach the design of $P$ is formulated as a hypothesis test with the null hypothesis $\mathcal{H}_0$ and the alternative hypothesis $\mathcal{H}_1$ \cite{willsky2003generalized}:
\begin{equation*}
   \begin{cases}
   \mathcal{H}_0: r(k) \sim \mathcal{N}(0,P\bar\Sigma_{0} P^\top)\\
       \mathcal{H}_1: r(k) \sim \mathcal{N}(PW_a a_s, P \bar \Sigma_{0}P^\top)
   \end{cases}  
\end{equation*}
where $\bar \Sigma_0 =  W_d \Sigma_0 W_d ^\top$. 
Given a preset FAR upper bound $\varepsilon$, GLRT gives rise to the threshold $J_{\rm th}=1$ and residual generator \eqref{eq_fd_parity} with the optimal design matrix
\begin{equation}
    P=\bar\Sigma_0^{-\frac{1}{2}}W_a(W_a^\top \bar\Sigma_0^{-1} W_a)^\dagger W_a^\top\bar\Sigma_0^{-{1}}/\sqrt{\chi_{m}^2(1-\varepsilon)}. \label{eq_GLRT}
\end{equation} 

\subsection{Wasserstein-based Distributionally Robust Detection Design }
To effectively detect cyber-attacks, a critical issue is how to describe the statistical properties of underlying disturbances and analyze its impact on detection performance. Under the Gaussian assumption, the well-established GLRT has been widely used to design cyber-attack detectors. However, the true disturbance distribution $\mathbb{P}_\xi$ is typically unknown and shows complicated characteristics, e.g. nonstationarity, heavy tails or multimodality, which deviate remarkably from Gaussianity. When such a distributional mismatch occurs, alarm floods can be induced, which makes practitioners eventually discredit the attack detector. 

To address uncertainty in probability distributions, we draw ideas from DRO to formulate the detector design problem, which captures the unknown disturbance distribution $\mathbb{P}_\xi$ by constructing a so-called ambiguity set $\mathcal{D}$ instead of a single parametric distribution. Specifically, as long as historical data are available for defining an empirical distribution, one can construct $\mathcal{D}$ as a family of probability distributions that are close to the empirical distribution at hand. To evaluate the distance between two probability distributions, the Wasserstein distance has been a popular option with clear interpretability and favorable statistical properties \cite{mohajerin2018data,shang2021distributionally}.  

\begin{definition}[Wasserstein distance, \cite{kantorovich1958space}]
For given two distributions $\mathbb{P},\mathbb{P}'\in \mathcal{M}(\Xi)$, where $\mathcal{M}(\Xi)$ is the probability space, whose random $\xi$ is supported on $\Xi$
\begin{equation*}
d_{\rm W}(\mathbb{P},{\mathbb{P}}') = \inf_{\mathbb{Q}\in\mathcal{Q}(\mathbb{P},{\mathbb{P}}')}  \mathbb{E}_{\mathbb{Q}}
\left\{ \lVert\xi-\xi'\rVert \right\},
\end{equation*}
where $\mathcal{Q}(\mathbb{P},{\mathbb{P}}') \in  \mathcal{M}(\Xi^2)$ is the set composed of all the joint distributions of $\xi$ and $\xi'$ with marginal distributions $\mathbb{P}$ and $\mathbb{P}'$.
\end{definition}

\begin{definition}[Wasserstein ambiguity set, \cite{mohajerin2018data}] \label{def_ambiguity}
Given $N$ independent samples $\{\hat{\xi}_i\}_{i=1}^N$, the Wasserstein ambiguity set is defined as:
\begin{equation*}
    \mathcal{D}_{\rm W}(\theta ; N) = \left\{\mathbb{P}\in \mathcal{M}(\Xi)\left|d_{\rm W}(\mathbb{P},\hat{\mathbb{P}}_{N})\leq\theta\right.\right\},
\end{equation*}
where $\mathcal{M}(\Xi)$ is the probability space supported on $\Xi$, $\theta$ is the radius of the ambiguity set, and $\hat{\mathbb{P}}_{N} = \frac{1}{N} \sum_{i=1}^{N} \delta_{\hat{\xi}_i}$ is the empirical distribution.
\end{definition}
Departing from the usual Gaussian assumption, the Wasserstein ambiguity set offers a general data-driven tool for characterizing uncertainty in the true distribution $\mathbb{P}_\xi$ provided that the radius $\theta$ is suitably chosen. Following the spirit of DRO, a Wasserstein distributionally robust (WDR) anomaly detector design problem can be formulated as follows \cite{shang2021distributionally}:
\begin{subequations}
\begin{align}
  \max_{P\succ0} &~ \lVert PW_a \rVert^2_F \label{eq_DRFD1}\\ 
   {\rm s.t.}&~  \sup_{\mathbb{P}_{\xi} \in \mathcal{D}_{\rm W}(\theta ; N)}\mathbb{P}_{\xi} \{ \xi^\top  W_d ^\top \bar P   W_d  \xi > 1 \} \leq \varepsilon \label{eq_DRFD2}
   \end{align}\label{eq_DRFD}
\end{subequations}
\noindent where ${\bar P = P^\top P\succ0}$ is defined to make the objective linear in $\bar P$. It maximizes the overall detectability $ \lVert PW_a \rVert^2_F$ in \eqref{eq_DRFD1} while respecting the distributionally robust chance constraint (DRCC) on FAR in \eqref{eq_DRFD2}. 

\begin{remark}
In engineering practice, non-stationarity often arises from factors such as aging-induced biases or gradual environmental shifts, manifesting as a continuous drift in the statistical properties of the disturbances. The FAR guarantee \eqref{eq_DRFD2} of the detector derived from \eqref{eq_DRFD} remains valid provided that the time-varying distribution $\mathbb{P}_\xi$ resides in the ambiguity set $\mathcal{D}_{\rm W}(\theta ; N)$. In fact, the observed violation of the FAR tolerance $\varepsilon$ can be seen as an indicator that $\mathcal{D}_{\rm W}(\theta ; N)$ cannot capture $\mathbb{P}_{\xi}$  due to the distributional non-stationarity, thereby calling for recalibration the Wasserstein radius $\theta$.
\end{remark}

\section{Main Result} \label{mainr}
\subsection{Stealthy Attack Set and Performance Degradation Metric}
The distributionally robust detector \eqref{eq_DRFD} only focuses on maximizing the overall detectability, while ignoring the adversarial nature of cyber-attacks. The attacker aims to cause damage to the CPS while avoiding raising security alarms. In this work, the attacker is assumed to have access to full system knowledge. Admittedly, for any detector design matrix $P$, there always exist stealthy attacks that are too unobvious to be detected. This leads to a high security risk that a well-designed stealthy attack may not trigger alarms but do harm to the CPS \cite{wu2022covert,zhang2022design}. Thus we seek to design a detector oriented towards minimizing performance degradation of the CPS subject to FDI attacks. 
Based on the residual generator and the threshold $J_{\rm th}=1$, we first define the stealthy attack set $\mathcal S_a$ as:
\begin{equation} 
\mathcal S_a=   \left\{a_s \left| \lVert P(W_d\xi+W_aa_s  )  \rVert ^2 \leq 1\right.\right\}. \label{eq_stealthy_set1}
\end{equation}
We assume that the attack projection matrix $W_a$ has full row rank, such that $ W_a W_a^\dagger=I$, and define a new attack representation $\bar a_s$: 
\begin{equation}
    a_s= -W_a^\dagger W_d\xi+ \bar a_s. \label{eq_bar_a}
\end{equation}
By substituting \eqref{eq_bar_a} into \eqref{eq_stealthy_set1}, one can decouple $\mathcal S_a$ from $\xi$, leading to the residual dynamic $r=PW_a \bar a_s$. On this basis, we can define an equivalent, but now deterministic stealthy attack set $\mathcal S_{\bar a}$ using  $\bar a_s$:
\begin{equation*} 
\mathcal S_{\bar a}=\left\{\bar a_s \left| \lVert PW_a \bar a_s\rVert ^2 \leq 1\right.\right\}.
\end{equation*}
However, the primary consequence of $W_a$ being row-rank-deficient is the existence of a non-trivial kernel space $\ker(W_a)$. Under this case, attackers can design strictly stealthy attacks $ a_s\in {\rm ker}(W_a)$, which satisfy $W_a a_s=0$ and thus will not be indicated at all by $r$ as per \eqref{eq_fd_parity}. Therefore, we decompose $\bar a_s$ into detectable components $\bar a^{\rm img}_s$  and strictly stealthy ones $\bar a^{\rm ker}_s$: 
\begin{align}
      & \bar a_s= \bar V_a\bar a^{\rm img}_s +\bar V_a^\perp
  \bar a^{\rm ker}_s, \label{eq_a_ker_img} \\
&   \bar V_a= V_{a,[1:sn_y-n_x]},~ \bar V^\perp_a =  V_{a,[sn_y-n_x+1:sn_a]} \nonumber
\end{align}
where the orthogonal matrix $V_a\in \mathbb{R}^{sn_a\times sn_a}$ comes from the singular value decomposition $W_a = U_a \Lambda V_a^\top$. Note that $\bar V_a^\perp \bar a^{\rm ker}_s$ is completely undetectable by \eqref{eq_fd_parity}, so we set $\bar a^{\rm ker}_s = 0$ and merely focus on the detectable component $\bar a^{\rm img}_s$. 
Consequently, the core of the stealthiness analysis reduces to characterizing the set of attacks $\bar a^{\rm img}_s$. 
By substituting \eqref{eq_a_ker_img} into \eqref{eq_fd_parity}, we arrive at $ r= PW_a\bar V_a\bar a^{\rm img}_s$, giving rise to the following stealthy attack set under study:
\begin{equation*} 
\mathcal S^{\rm img}_{\bar a}=\left\{\bar a^{\rm img}_s \left| \lVert  PW_a\bar V_a\bar a^{\rm img}_s\rVert ^2 \leq 1\right.\right\}.
\end{equation*}

Next we analyze the effect of stealthy attacks on system states in the absence of the true distribution of $\xi$. The augmented state $\bar{x}$ in \eqref{eq_SS2} can be decomposed as follows: 
\begin{align}
\begin{cases}
   \bar {x}(k)=\bar {x}^{\rm nom}(k)+\bar {x}^{\rm dev}(k)\\
\bar {x}^{\rm nom}(k+1)=\bar A{\bar x^{\rm nom}}(k)+ \bar B_{ r} {y_{\rm ref}(k)} \\
\bar {x}^{\rm dev}(k+1)=\bar A{\bar x^{\rm dev}}(k)+\bar B_a a(k)+\bar B_d d(k)
\end{cases} \label{eq_x2}
\end{align}
where $\bar {x}^{\rm nom}$ represents the deterministic response to the given reference trajectory $y_{\rm ref}$, and $\bar {x}^{\rm dev}$ captures the state deviation caused by both the attack $a$ and uncertain disturbance $d$. The reachable set of $\bar {x}^{\rm dev}$ can be useful for evaluating performance degradation under stealthy attack. To align with the $s$-long stealthy attacks $\bar a_s$ in $\mathcal S_{\bar a}$, we first augment the dynamics \eqref{eq_x2} of the deviation component over a time horizon of $s$: 
\begin{equation} {\bar x^{\rm dev}_{s}}(k+s)=\bar A_s{\bar x^{\rm dev}_{s}}(k)+\bar B_{ a,s}  a_s(k+s-1)+\bar B_{d,s} \xi(k+s-1) \label{eq_SS_dev}
\end{equation}
where 
\begin{align*}
       &\bar A_s =\begin{bmatrix}
        0&0&\cdots&  \bar{A}\\
                0&0&\cdots&  \bar{A}^2\\
       \vdots& \vdots&\ddots &\vdots\\
      0&0&\cdots& \bar{A}^s
    \end{bmatrix},\\
    & \bar B_{ a,s}=\begin{bmatrix}
    \bar B_{ a}&0&\cdots&0\\
    \bar A_s\bar B_{a}&\bar B_{a}&\cdots&0\\
    \vdots&\vdots&\ddots&\vdots\\
    \bar A_s^{s-1}\bar B_{a}&  \bar A_s^{s-2}\bar B_{a}&\cdots&\bar B_{a}\end{bmatrix},
\end{align*}
and $\bar B_{d,s}$ are defined akin to $\bar B_{ a,s}$ with $\bar B_a$ replaced by $\bar B_d$. One can further replace $a_s$ with $\bar a_s$ by substituting \eqref{eq_bar_a} and \eqref{eq_a_ker_img} into the augmented system \eqref{eq_SS_dev}:
\begin{equation}
\begin{aligned}
    {\bar x^{\rm dev}_{s}}(k+s)&=\bar A_s{\bar x^{\rm dev}_{s}}(k)+\bar B^{\rm img}_{\bar a,s} \bar a^{\rm img}_s(k+s-1)\\
    &+\bar B^{\rm ker}_{\bar a,s} \bar a^{\rm ker}_s(k+s-1)+\bar B_{\xi} \xi(k+s-1) 
\end{aligned}
\label{eq_SS_dev2}
\end{equation}
where
\begin{align*}
    &\bar B^{\rm img}_{\bar a,s}=\bar B_{ a,s} \bar V_a,~
    \bar B^{\rm ker}_{\bar a,s}=\bar B_{ a,s} \bar V_a^\perp\\
    &\bar B_{\xi} = \bar B_{ d,s}-\bar B_{ a,s}W_a^\dagger W_d.
\end{align*}

Since the distributions $\mathbb{P}_\xi$ in the ambiguity $\mathcal{D}_{\rm W}(\theta ; N)$ may have an unbounded support, i.e., $\Xi=\mathbb{R}^{n_\xi}$, traditional deterministic set-based reachability analysis is not directly applicable. Consequently, we turn to a distributionally robust confidence disturbance set $\mathcal{E}_\xi(\beta) = \left\{ \xi \left| \lVert \xi\rVert^2_{Q(\beta)}\leq 1\right.\right\}$ with $Q(\beta)\succ0$ and a high confidence level $\beta \in (0,1)$ (e.g., $0.95$), satisfying the following DRCC:
\begin{equation}
\begin{aligned}
&\inf_{\mathbb{P}_{\xi} \in \mathcal{D}_{\rm W}(\theta, N)}  \mathbb{P}_{\xi}\left\{\xi \in \mathcal{E}_\xi(\beta)\right\}\\
= &\inf_{\mathbb{P}_{\xi} \in \mathcal{D}_{\rm W}(\theta, N)}  \mathbb{P}_{\xi}\left\{ \lVert \xi\rVert^2_{Q(\beta)} \leq 1\right\}\geq \beta.
\end{aligned} \label{eq_DRCC_Q_beta}
\end{equation}
The DRCC \eqref{eq_DRCC_Q_beta} ensures that $\mathcal{E}_\xi(\beta)$ is always a valid bounded $\beta$-confidence set for any distribution $\mathbb{P}_\xi$ residing in the Wasserstein ball $\mathcal{D}_{\rm W}(\theta,N)$. Note that rather than a predefined constant matrix, $Q$ is a decision variable to be optimized over in our design problem, and thus implicitly depends on $\beta$. So here we make clear the dependence of $Q(\beta)$ on beta, which will be omitted for brevity if no confusion is made. Based on the bounded set $\mathcal{E}_\xi(\beta)$, we define the following tightened deterministic reachable set at time $k$ under stealthy attack: 
\begin{equation*}
\begin{aligned}
    & \mathcal{R}_ {\bar x^{\rm dev}_{s}}(k)=
     \left\{  \bar x^{\rm dev}_{s}(k)~\middle|~
    \begin{aligned}
        &  \eqref{eq_SS_dev2}, ~\bar a^{\rm ker}_s(i) = 0, ~  i \in \mathbb{N}_{1:k}, \\
    &  \lVert  PW_a\bar V_a\bar a^{\rm img}_s(i)\rVert ^2 \leq 1 , ~ i \in \mathbb{N}_{1:k}, \\
       &\lVert \xi(i)\rVert^2_{Q}\leq 1, ~ i \in \mathbb{N}_{1:k}.
    \end{aligned}
    \right\},
    \end{aligned}
\end{equation*}
which characterizes possible state deviations caused by both stealthy attacks and disturbance in the $\beta$-confidence set. Intuitively, a larger size of $\mathcal{R}_ {\bar x^{\rm dev}_{s}}(k)$ indicates severer performance degradation and consequently, a heavier impact of stealthy attacks on CPS security; however, its explicit dependence on $k$ erects obstacles for evaluating the performance degradation. Therefore, we consider its \textit{asymptotic outer approximation} in the form of an ellipsoid irrespective of $k$:
\begin{equation}
 \mathcal{R}_ {\bar x^{\rm dev}_{s}}(k)
 \subseteq \mathcal{E}^{\infty}_ {\bar x^{\rm dev}_{s}}
 = \left\{  \bar x^{\rm dev}_{s}~\middle| \lVert \bar x^{\rm dev}_{s} \rVert_M ^2 \leq 1 \right\},~ k  \rightarrow \infty \label{eq_outer_appro}
\end{equation}
where $M\succ0$ decides both the shape and the volume of the ellipsoid $\mathcal{E}^{\infty}_ {\bar x^{\rm dev}_{s}}$. In a nutshell, $\mathcal{E}^{\infty}_ {\bar x^{\rm dev}_{s}}$ provides an approximation of the asymptotical reachable set with stochastic disturbances from the $\beta$-confidence set $\mathcal{E}_\xi(\beta)$. It is worth noting that a smaller volume of $\mathcal{E}_\xi(\beta)$ means tighter bounds on state deviations, implying better containment of attack effects. It is known that $1/\sqrt{\det(M)}$ is proportional to the volume of $\mathcal{E}^{\infty}_ {\bar x^{\rm dev}_{s}}$ \cite{boyd1994linear}, which naturally motivates the utilization of $-\log \det(M)$ as a metric for quantifying performance degradation under stealthy attacks and eventually guiding the design of $P$. 

\subsection{Derivation of Asymptotical Outer Approximation}
Our objective is to find the shape matrix $M$ that defines the smallest possible invariant ellipsoid $\mathcal{E}^{\infty}_ {\bar x^{\rm dev}_{s}}$. To this aim, we first present a preliminary result, which is useful for developing an asymptotic bound on the state evolution under multiple bounded inputs.
\begin{lemma} \label{lemma_set}(\cite[Lemma 1]{murguia2020security}).
Suppose a given constant $\alpha \in (0, 1)$ and a nonnegative function $V(k)$, if $\omega_i(k)^\top M_i\omega_i(k) \leq 1$ with $M_i \succ 0$, $i\in \mathbb{N}_{1:N_{\alpha}}$ and there exist $\alpha_i\in (0, 1),~i\in \mathbb{N}_{1:N}$, satisfying $\sum_{i=1}^{N_{\alpha}} \alpha_i \geq \alpha$ and
\begin{equation}
      V({k+1}) \leq \alpha V(k) + \sum_{i=1}^{N_{\alpha}} (1 - \alpha_i)\omega^\top_i(k) M_i(k) \omega_i(k), \label{eq_V_set}
\end{equation}
then, the bound of  $V(k)$ satisfies
\begin{equation*}
V({k}) \leq  \alpha^{k-1} V(1) + \frac{(N_{\alpha} - \alpha)(1 - \alpha^{k-1})}{1 - \alpha},
\end{equation*}
and its asymptotic upper bound is given by $\lim_{k \to \infty} V(k) \leq {(N_{\alpha} - \alpha)}/{(1 - \alpha)}$.
\end{lemma}
Now we are in a position to derive the necessary conditions that the matrix $ M$ must satisfy to ensure the ellipsoid $\mathcal{E}^{\infty}_ {\bar x^{\rm dev}_{s}}$ covers. Therefore, Theorem \ref{prop1} is proposed to formalize these conditions.

\begin{theorem} \label{prop1}
  For a given $\alpha \in (0,1)$ if there exist $\alpha_1$, $\alpha_2$ and $\bar M$ satisfying
\begin{equation} 
\begin{cases}
&\bar M = \frac{2-\alpha}{1-\alpha} M \succ 0,~ \alpha_1+\alpha_2\geq \alpha\\
      & 0 \leq \alpha_1 \leq 1, 0 \leq \alpha_2  \leq 1\\ 
  &  \begin{bmatrix}
        \alpha\bar M&0&\bar A_s^\top \bar M\\  
        0& \bar \Omega&\bar B_s^\top \bar M \\
        \bar M \bar A_s&\bar M \bar B_s &\bar M
    \end{bmatrix} \succeq 0 
\end{cases} \label{eq_prop1}
    \end{equation}
   where 
\begin{align*}
   &  \bar\Omega = {\rm diag}\left\{ (1-\alpha_1) \bar V_a^\top W_a^\top \bar{P}W_a  \bar  V_a,(1-\alpha_2)Q\right\}, \\
& \bar B_s = \begin{bmatrix}
\bar B^{\rm img}_{\bar a,s}&\bar B_{\xi}
\end{bmatrix},
\end{align*}
\eqref{eq_outer_appro} with $M$ is a valid asymptotic outer approximation of the reachable set $\mathcal{R}_ {\bar x^{\rm dev}_{s}}(k)$.

\begin{proof}
We first define $\zeta(k)=\bar x^{{\rm dev}}_{s} (sk)$, 
$V(k) = \zeta^\top(k) \bar M \zeta(k) $,  $\omega_1(k)= \bar a_s^{\rm img}(s(k+1)-1)$, $\omega_2(k)=\xi(s(k+1)-1)$, $\omega(k)=\begin{bmatrix}
    \omega_1(k)   &  \omega_2(k)
\end{bmatrix}$ and $N_{\alpha}=2$. If there exists $\alpha_1+\alpha_2 \geq \alpha \in (0,1)$, such that \eqref{eq_V_set} holds for any $\bar x^{{\rm dev}}_{s} (sk)$. By substituting \eqref{eq_SS_dev2} into \eqref{eq_V_set}, we obtain 
\begin{equation*}
\begin{aligned}
& \alpha \zeta^\top (k)\bar M  \zeta (k) -
\big(\bar A_s \zeta(k)+\bar B_{s}  {\omega}(k) \big)^\top \bar M \big(\bar A_s \zeta(k)+\bar B_{s}  {\omega}(k) \big)\\
&+ (1-\alpha_1) \omega_1^\top(k)  W_a^\top \bar{P}W_a  \omega_1(k)+(1-\alpha_2) \omega_2^\top(k)Q\omega_2(k)
  \\= & \begin{bmatrix}
       \zeta(k) \\ \omega(k)
\end{bmatrix}^\top \begin{bmatrix}
        \alpha\bar M-\bar A_s^\top \bar M \bar A_s&- \bar A_s^\top \bar M \bar B_{s}\\  
      - \bar B_s^\top \bar M \bar A_{s}  &\bar \Omega-B^\top_{s} \bar M \bar B_{s}
    \end{bmatrix} 
\begin{bmatrix}
       \zeta(k) \\\omega(k)
\end{bmatrix}\geq0,
\end{aligned}
\end{equation*}
which implies
\begin{equation}
   \begin{bmatrix}
        \alpha\bar M-\bar A_s^\top \bar M \bar A_s&- \bar A_s^\top \bar M \bar B_{s}\\  
      - \bar B_s^\top \bar M \bar A_{s}  &\bar \Omega-B^\top_{s} \bar M\bar B_{s}
    \end{bmatrix} \succeq 0. \label{eq_schur}
\end{equation}
Due to $M \succ 0$, we can apply the Schur complement to \eqref{eq_schur} to derive \eqref{eq_prop1}. Based on Lemma \ref{lemma_set}, this establishes that for the state sequence $\bar{x}_s^{\text{dev}}(sk)$, the asymptotic upper bound $\lim_{k\to\infty} V(k) \leq \frac{2-\alpha}{1-\alpha}$ holds. Due to $\alpha \in (0,1)$, the influence of any initial states $\{\bar x^{{\rm dev}}_{s} (i) \}^{s-1}_{i=0}$ vanishes exponentially over time. Therefore, the asymptotical outer approximation of the reachable set $\left\{ \bar{x}^{\rm dev}\left |\bar{x}^{{\rm dev}\top}M \bar{x}^{\rm dev}\leq1\right.\right\}$ can be extended to the entire trajectory $\bar{x}^{\rm dev}(k)$, thus completing the proof.
\end{proof}
\end{theorem}

\subsection{Design Problem Formulation and Solution Algorithm}
The proposed distributionally robust design of FDI attack detector is formalized as follows: 
\begin{subequations}\label{eq_DRAD_Primal}
\begin{align} 
\min_{\begin{subarray}{c}
\bar P \succ 0,Q\succ 0,\\\bar M\succ0,\alpha_1,\alpha_2
\end{subarray}} &~ -{\log\det}( \bar M)  \label{eq_design1}\\
{\rm s.t.}~~~ &~  \sup_{\mathbb{P}_{\xi} \in \mathcal{D}_{\rm W}(\theta, N)}
\mathbb{P}_\xi \left\{\lVert W_d \xi\rVert_{\bar P}^2> 1 \right\} \leq \varepsilon \label{eq_design2}\\
&~ \sup_{\mathbb{P}_{\xi} \in \mathcal{D}_{\rm W}(\theta, N)}
\mathbb{P}_\xi \left\{\lVert \xi\rVert_{Q}^2> 1 \right\} \leq1- \beta \label{eq_design3}\\
&~ {\rm Constraint}~ \eqref{eq_prop1} \label{eq_design4}
\end{align} \label{eq_design}
\end{subequations}
In this design,  the FAR is guaranteed to be below the tolerance level $\varepsilon $ for any disturbance distribution $\mathbb{P}_\xi$ within the ambiguity set $\mathcal{D}_{\rm W}$ by the DRCC \eqref{eq_design2}. The DRCC \eqref{eq_design3} defines a safe $\beta$-confidence set for $\xi$ that is useful for deriving the asymptotic outer approximation in \eqref{eq_design4}. Constraint \eqref{eq_design4} derived in Theorem \ref{prop1} is the core of defining the proposed security metric by finding an ellipsoid parameterized by $\bar M$ being a valid outer approximation of the asymptotic reachability set of state deviations. Note that minimizing the proposed metric $-\log \det(M)$ is equivalent to minimizing $-{\log\det}( \bar M )$. Thus we adopt the convex objective function $-{\rm logdet}( \bar M )$ as the objective, which seeks to reduce the CPS state deviation caused by stealthy attacks as much as possible. However, the optimization problems \eqref{eq_design2} and \eqref{eq_design3} contain semi-infinite optimization problems, which are computationally intractable. To handle two Wasserstein-based DRCCs, we introduce the following result.

\begin{lemma}\label{lemma_DRCC}(\cite[Theorem 5]{shang2021distributionally}).
For a given probability threshold $\varepsilon$ and $\Xi =\mathbb{R}^{n_\xi} $, the worst-case FAR chance-constraint over the Wasserstein-based ambiguity set $\mathcal{D}_{\rm W}(\theta, N)$
\begin{equation*}
    \sup_{\mathbb{P}_{\xi} \in \mathcal{D}_{\rm W}(\theta, N)}  \mathbb{P}\left\{\lVert W_d\xi\rVert^2_{\bar P} > 1\right\}\leq \varepsilon\\
\end{equation*}
holds if and only if 
\begin{align}
\begin{cases}
    \displaystyle& \lambda \geq 0 \\
 \displaystyle&  y_i\geq 0,~\tau_i\geq 0,~ y_i\geq\lambda-t_i,~ i \in \mathbb{N}_{1:N}\\
    \displaystyle&\theta + \frac{1}{N}\sum_{i=1}^N y_i \leq \lambda\varepsilon,~  i \in \mathbb{N}_{1:N}\\
 \displaystyle & \begin{bmatrix} 
   I  & -\hat{\xi}_i \\
   -\hat{\xi}_i^\top & \hat{\xi}_i^\top \hat{\xi}_i - q_i 
   \end{bmatrix} -\tau_i
   \begin{bmatrix}
         W_d^\top \bar{P} W_d &0\\0&-1
   \end{bmatrix}\succeq 0,\\
   & \qquad\qquad\qquad\qquad\qquad\qquad\qquad\quad ~~ i \in \mathbb{N}_{1:N}, \\
   \displaystyle &\begin{bmatrix} 
   q_i & t_i \\
   t_i & 1 
   \end{bmatrix} \succeq 0,~  i \in \mathbb{N}_{1:N}.
\end{cases}\label{eq_lemma_DRCC}
\end{align}

\end{lemma}
In the light of Lemma \ref{lemma_DRCC}, we can now address the semi-infinite constraints \eqref{eq_design2} and \eqref{eq_design3} in our design problem. Therefore, Theorem \eqref{theo_reformu} is proposed by replacing them with a finite set of constraints. 
\begin{theorem} \label{theo_reformu}
The distributionally robust attack detector design problem \eqref{eq_design} with $\Xi=\mathbb{R}^{n_{\xi}}$ admits the following exact reformulation: 
\begin{equation}
    \begin{aligned}
        \min_{\begin{subarray}{c}\bar{P},\bar{M},{Q},\alpha_1,\alpha_2\\
        \gamma,\lambda,y_i,v_i,\tau_i,\\
        \pi_i,t_i,u_i,p_i,q_i\end{subarray}}  &~ -{ \log\det}( \bar M)  \\
{\rm s.t.}~~~~ &~  \gamma \geq 0,~u_i\geq 0,~\pi_i\geq 0, ~ i \in \mathbb{N}_{1:N}\\
    &~\theta + \frac{1}{N}\sum_{i=1}^N u_i \leq(1-\beta)\gamma\\
&~ \begin{bmatrix} 
   I & -\hat{\xi}_i \\
   -\hat{\xi}_i^\top & \hat{\xi}_i^\top \hat{\xi}_i - p_i 
   \end{bmatrix} -\pi_i \begin{bmatrix}
       Q&0\\0&-1
   \end{bmatrix} \succeq 0,\\
   ~& \qquad\qquad\qquad\qquad\qquad \qquad   \quad  i \in \mathbb{N}_{1:N} \\
   &~\begin{bmatrix} 
   p_i & u_i \\
   u_i & 1 
   \end{bmatrix} \succeq 0, i \in \mathbb{N}_{1:N}\\
&~ \bar P\succ 0,~Q\succ 0\\
&~ {\rm Constraints}~\eqref{eq_prop1},~ \eqref{eq_lemma_DRCC}
    \end{aligned}\label{eq_reformu}
\end{equation}
\begin{proof}
The proof is a direct application of Theorem \ref{prop1} and Lemma \ref{lemma_DRCC}, so it is omitted for brevity.
\end{proof}
\end{theorem}
A key observation is that if either set of bilinear variables $\{\tau_i,\pi_i,a_1,a_2\}$ is fixed, the problem becomes a semi-definite program (SDP), which is convex with respect to the other set $\{\bar P,Q\}$ and thus becomes solvable. To address the NP-hard problem (22), we employ the sequential minimization based on coordinate descent, which is a standard approach to solving BMIs \cite{el1994synthesis,hassibi1999path}. Although this strategy ensures a non-decreasing sequence of objective values, it does not guarantee global convergence, yet it has the potential to yield a high-quality solution. 
To implement the sequential minimization, we first rewrite \eqref{eq_reformu} as follows:
\begin{equation}
\begin{split} 
\min_{\bar{P} \succ 0,  {Q} \succ 0} &~ -{\log\det}(\bar M) \\
{\rm s.t.}~~&~ \mathcal{J}_1(\bar{P}) \leq 0,~ \mathcal{J}_2(Q) \leq 0,\\
&~ \mathcal{J}_3(\bar P,Q) \geq a
\end{split}\label{eq_J0}
\end{equation}
where $\mathcal{J}_1(\bar{P})$ is the optimal value of the following subproblem:
\begin{equation}
\begin{split}
\min_{\lambda, y_i, \tau_i, t_i, q_i}&~ \theta + \frac{1}{N} \sum^{N}_{i=1}y_i - \lambda \varepsilon \\
{\rm s.t.}~~~&~y_i \geq 0 ,~y_i \geq \lambda  - t_i,~\tau_i\geq0,~ i\in\mathbb{N}_{1:N}\\
& \begin{bmatrix} 
   I - \tau_i W_d^\top \bar{P} W_d & -\hat{\xi}_i \\
   -\hat{\xi}_i^\top & \hat{\xi}_i^\top \hat{\xi}_i - q_i + \tau_i 
   \end{bmatrix} \succeq 0,\\
&~ \qquad \qquad \qquad \qquad \qquad \qquad  \qquad  i\in\mathbb{N}_{1:N} \\
&~\begin{bmatrix}q_i&t_i\\t_i&1\end{bmatrix}\succeq 0,~ i\in\mathbb{N}_{1:N},~\lambda \geq 0
\end{split}\label{eq_J1}
\end{equation}
$ \mathcal{J}_2(Q)$ is the optimal value of the following subproblem:
\begin{equation}
\begin{split}
\min_{\gamma,v_i, \pi, u_i, p_i}&~ \theta + \frac{1}{N}\sum_{i=1}^N v_i -(1-\beta)\gamma \\
{\rm s.t.}~~~&~v_i \geq 0 ,~v_i \geq \eta -  u_i ,~\pi_i\geq0,~ i\in\mathbb{N}_{1:N}\\
&~ \begin{bmatrix} 
   I-\pi_iQ & -\hat{\xi}_i \\
   -\hat{\xi}_i^\top & \hat{\xi}_i^\top \hat{\xi}_i - p_i +\pi_i
   \end{bmatrix}  \succeq 0,~ i\in\mathbb{N}_{1:N} \\
 &~\begin{bmatrix} 
   p_i & u_i \\
   u_i & 1 
   \end{bmatrix} \succeq 0,~ i\in\mathbb{N}_{1:N},~\gamma \geq 0
   \end{split}\label{eq_J2}
\end{equation}
and $ \mathcal{J}_3(\bar P,Q)$ is the optimal value of the following subproblem: 
\begin{equation}
\begin{split}
\max_{\alpha_1,\alpha_2,\bar M}&~ \alpha_1+\alpha_2\\
{\rm s.t.}~~&~  0 \leq \alpha_1 \leq 1,~0 \leq \alpha_2 \leq 1,~ \bar M \succ 0\\
  &~   \begin{bmatrix}
    \alpha     \bar M&0&\bar A_s^\top \bar M\\  
        0& \bar \Omega&\bar B_s^\top \bar M \\
        \bar M \bar A_s&\bar M \bar B_s &\bar M
    \end{bmatrix} \succeq 0 \\
   & ~ \bar\Omega = {\rm diag} \left \{
 (1-\alpha_1) \bar V_a^\top W_a^\top \bar{P}W_a  \bar  V_a,(1-\alpha_2)Q \right \}
\end{split}\label{eq_J3}
\end{equation}
This decomposition allows us to perform sequential minimization to resolve the original non-convex problem \eqref{eq_reformu}. First, we initialize $\bar P$ and $Q$ with small positive-definite matrices to ensure initial feasibility. 
Then, we can alternately solve three SDPs \eqref{eq_J1}, \eqref{eq_J2} and \eqref{eq_J3} with fixed $\{\bar P,Q\}$ to enlarge the feasible region, and solve the SDP \eqref{eq_J0} with fixed $\{\tau_i,\pi_i,\alpha_1,\alpha_2\}$ to achieve a higher objective. This procedure is repeated until convergence or the maximum number of iterations. Notably, these SDPs can be solved using off-the-shelf solvers. The worst-case computational complexities of \eqref{eq_J0}, \eqref{eq_J1}, \eqref{eq_J2} and \eqref{eq_J3} are $\mathcal{O}\left(\left[s(n_y+n_x+n_{\xi})\right]^{6}\right)$, $\mathcal{O}\left(\left(sn_y\right)^{6}\right)$, $\mathcal{O}\left(\left(sn_{\xi}\right)^{6}\right)$ and $\mathcal{O}\left(\left(sn_x\right)^{6}\right)$ when $N=\mathcal{O}\left(s\left(n_y+n_x+n_{\xi} \right)\right)$.  Nevertheless, the proposed method is an offline design procedure and does not require real-time execution.
As this iterative scheme is a form of coordinate descent, the objective function is guaranteed to be non-decreasing at each iteration and finally converges. The full implementation details are summarized in Algorithm \ref{alg1}. 

\begin{algorithm}
\caption{Solution algorithm for reformulation of the distributionally robust detector design problem \eqref{eq_reformu}}
\label{alg1}
\begin{algorithmic}[1]
\Require Coefficient matrices $W_a$, $W_d$, $\bar A_s$, $\bar B_{\bar a,s}$, $\bar B_{\xi}$, disturbance samples $\{\hat{\xi}_i\}_{i=1}^N$, Wasserstein radius $\theta >0$, tolerable FAR $0<\varepsilon<1$ and asymptotic convergence rate $0<\alpha<1$.  
\State Initialize $\bar{P}^{(0)} \leftarrow  \epsilon I$ and $Q \leftarrow  \epsilon I$ with a small enough $\epsilon>0$ and $n_{\rm ite}\leftarrow 0$. 
\While{stopping criteria not met}
\State Solve \eqref{eq_J1} over $\{\lambda, y_i, \tau_i, t_i, q_i\}$ with $\bar{P} \leftarrow \bar{P}^{(n_{\rm ite})}$, and obtain the optimal solution $\tau_i^{(n_{\rm ite})}$.
\State Solve \eqref{eq_J2} over $\{\gamma,v_i, \pi, u_i, p_i\}$ with $Q  \leftarrow Q^{(n_{\rm ite})}$, and obtain the optimal solution $\pi_i^{(n_{\rm ite})}$.
\State Solve \eqref{eq_J3} over $\{\alpha_1,\alpha_2,\bar M\}$ with $\bar{P}  \leftarrow \bar{P}^{(n_{\rm ite})}$ and $Q\leftarrow Q^{(n_{\rm ite})}$, and obtain the optimal solution $\alpha_1^{(n_{\rm ite})}$ and $\alpha_2^{(n_{\rm ite})}$. 
\State Solve \eqref{eq_J0} with $\tau_i \leftarrow \tau_i^{(n_{\rm ite})}$, $\pi_i \leftarrow \pi_i^{(n_{\rm ite})}$, $\alpha_1 \leftarrow \alpha_1^{(n_{\rm ite})}$ and $\alpha_2 \leftarrow \alpha_2^{(n_{\rm ite})}$, and obtain the optimal solution $\bar{{P}}^{(n_{\rm ite}+1)}$ and  $Q^{(n_{\rm ite}+1)}$ . 
 \State $n_{\rm ite} \leftarrow n_{\rm ite} + 1$.
\EndWhile
\State \textbf{Return} $\bar{P}$ .
\end{algorithmic}
\end{algorithm}
\begin{remark}
A practical simplification to the design problem \eqref{eq_DRAD_Primal} can be made by coupling the DRCC on FAR \eqref{eq_design2} and the DRCC on disturbance confidence set \eqref{eq_design2}. 
That is, letting $\beta=1-\varepsilon>0$, then we can choose the disturbance confidence set in asymptotic reachability analysis to be the same set of disturbances that do not trigger a false alarm, i.e., $\mathcal{E}_\xi(\beta) = \left\{ \xi \left| \lVert W_d\xi\rVert^2_{\bar P} \leq 1\right.\right\}$. This strategy can effectively eliminate the need to solve for the matrix $Q$ and the associated BMIs. Consequently, the subproblem \eqref{eq_J2} vanishes in each iteration, thereby improving the computational efficiency of  Algorithm \ref{alg1}.
\end{remark}

\begin{remark} \label{rmk_hyper}
The proposed design involves four key hyperparameters, $\{\varepsilon, \beta, \alpha,\theta\}$, whose selection requires practical consideration. First, the parameters $\varepsilon$ and $\beta$ possess clear physical interpretations, representing the FAR tolerance and the confidence level for the distributionally robust disturbance set $\mathcal{E}_{\xi}(\beta)$, respectively. 
This allows users to straightforward select their values based on their preference on practical tolerance of false alarms. Next, one can simply choose $\beta=1-\varepsilon$ for convenience. $\alpha$ is the asymptotic convergence rate in the derivation of $\mathcal{E}^{\infty}_ {\bar x^{\rm dev}_{s}}$. Because $\alpha \in  (0,1)$, its tuning is conducted via the well-established grid search method in \cite{zhang2024design,murguia2020security}. Then the Wasserstein radius theta can be effectively selected using $K$-fold cross-validation, which is a standard practice in Wasserstein DRO \cite{shang2021distributionally,feng2024distributionally}. Specifically, we define a search grid for $\theta$ and partition the empirical samples $\{\xi_i\}_{i=1}^N$ into $K$ folds. For each candidate $\theta$, we derive the design \eqref{eq_reformu} using Algorithm \ref{alg1} based on $K-1$ folds while testing the empirical FAR on the excluded fold. By repeating this procedure for each individual test fold, the optimal $\theta$ is decided as the smallest value such that the average FAR remains below the tolerance $\varepsilon$, which ensures a safe design with reduced conservatism.
\end{remark}

\subsection{Design under Row Rank-Deficient $W_a$}
In practice, the attack projection matrix $W_a$ may not always have full row rank due to the sensor redundancy or the inherent actuator structure. This section is therefore organized to address this more general and challenging case. We continue to employ the presentation \eqref{eq_bar_a} and \eqref{eq_a_ker_img}, but redefine 
\begin{align}
\bar V_a= V_{a,[1:{\rm rank}(W_a)]},~ \bar V^\perp_a =  V_{a,[{\rm rank}(W_a)+1:sn_a]} 
\end{align}
because $W_a$ is row rank-deficient. By substituting \eqref{eq_a_ker_img} into \eqref{eq_fd_parity}, the residual $r$ is governed by 
\begin{equation}
 r=  P(I-  W_aW_a^\dagger)W_d\xi+  PW_a\bar V_a\bar a^{\rm img}_s. \label{eq_res_a_img}
\end{equation}
Because of $  W_aW_a^\dagger \neq I $, \eqref{eq_res_a_img} only eliminates the component of $W_d\xi$ that lies in the row space of $W_a$, but it cannot completely remove the influence of $\xi$. To characterize a deterministic set of stealthy attacks 
$\bar a^{\rm img}_s$, we adopt two distinct design philosophies, leading to two different definitions. First, a conservative stealthy attack set $\bar{\mathcal S}^{\rm img}$ is defined as:
\begin{equation} 
\bar{\mathcal S}^{\rm img}_{\bar a,1}=\left\{\bar a_s^{\rm img} \left|
\begin{aligned}   
    \lVert(I-W_a W_a^\dagger) W_d\xi+  W_a\bar V_a\bar a^{\rm img}_s\rVert_{\bar P} ^2 \leq 1,\\
    \forall \xi \in \mathcal{E}_\xi(\eta)
\end{aligned}
\right.\right\}, \label{eq_set_forall}
\end{equation}
which only includes those attacks that can remain undetectable for every possible realization of disturbance $\xi$ in the confidence set $\mathcal{E}_\xi(\eta)$. From the defender's point of view, this is an optimistic assumption that the attacker will only launch attacks that are guaranteed to be stealthy with a probability of at least $\eta$. Next, we develop a more tractable approximation of \eqref{eq_set_forall}.

\begin{theorem} \label{theo_conser}
By choosing $\eta = 1-\varepsilon>0$ and $\mathcal{E}_\xi(\eta) = \left\{ \xi \left| \lVert W_d\xi\rVert^2_{\bar P}  \leq 1\right.\right\}$ being a valid $\eta$-confidence set, if
\begin{equation} 
\lVert W_a\bar V_a\bar a^{\rm img}_s\rVert_{\bar P} ^2 \leq \frac{1}{2}- \lVert I-W_a W_a^\dagger \rVert_2^2 \label{eq_strict2}
\end{equation}
then
\begin{equation}
\lVert (I-W_a W_a^\dagger) W_d\xi+  W_a\bar V_a\bar a^{\rm img}_s\rVert_{\bar P} ^2 \leq 1, ~\forall \xi \in \mathcal{E}_\xi(\eta).  \label{eq_strict1}
\end{equation} 
\begin{proof}
To derive a tractable sufficient condition, we upper bound for the  left hand of \eqref{eq_strict1}:
    \begin{equation}
    \begin{aligned}
  &~ \lVert(I-W_a W_a^\dagger) W_d\xi+  W_a\bar V_a\bar a^{\rm img}_s\rVert_{\bar P} ^2\\
\leq&~ 2  \lVert(I-W_a W_a^\dagger) W_d\xi\rVert_{\bar P} ^2+2
\lVert W_a\bar V_a\bar a^{\rm img}_s\rVert_{\bar P} ^2\\
\leq &~2  \lVert I-W_a W_a^\dagger\rVert_2^2 \lVert  W_d\xi\rVert_{\bar P} ^2+2
\lVert W_a\bar V_a\bar a^{\rm img}_s\rVert_{\bar P} ^2\\
\leq&~2  \lVert I-W_a W_a^\dagger\rVert_2^2 +2
\lVert W_a\bar V_a\bar a^{\rm img}_s\rVert_{\bar P} ^2\\
\end{aligned}\label{eq_strict3}
\end{equation}
Here, the first inequality follows from the Cauchy-Schwarz inequality, the second inequality is due to the property of the spectral norm of matrix, and the last inequality holds by the definition of  $\mathcal{E}_\xi(\eta)$. By inserting \eqref{eq_strict2} into \eqref{eq_strict3}, we obtain \eqref{eq_strict1}, which completes the proof.
\end{proof}
\end{theorem}
\begin{remark}
It is worth noting the tightness of this sufficient condition \eqref{eq_strict2}. 
A noteworthy special case occurs when the singular value of $ I-W_a W_a^\dagger $ is large enough to make the right-hand side of $\lVert I-W_a W_a^\dagger \rVert^2_2>\frac{1}{2}$, leading to an empty $\bar{\mathcal S}^{\rm img}_{\bar a}$. The physical interpretation is that there exists a large degree of freedom for the disturbances $\xi$ that are orthogonal to the row space $W_a$ to affect the residual $r$. Consequently, the worst-case disturbance becomes so influential that the alarm can be triggered under any attacks $a^{\rm img}_s$.
\end{remark}
To alleviate conservatism in ${\mathcal S}_{\bar a,1}^{\rm img}$, we introduce an alternative stealthy attack set ${\mathcal S}^{\rm img}_{\bar a,2}$ as:
\begin{equation*} 
{\mathcal S}^{\rm img}_{\bar a,2}=\left\{\bar a_s^{\rm img} \left|
\begin{aligned}
   \lVert W_a\bar V_a\bar a^{\rm img}_s\rVert_{\bar P} ^2 \leq 1
\end{aligned}
\right.\right\}, 
\end{equation*}
which is a broader and more computationally tractable set that characterizes all attacks that are stealthy in the disturbance-free case. This seemingly simple set also bears a clear probabilistic interpretation. If the disturbance $\xi$ is drawn from any centrally symmetric distribution, any attacks $\bar a^{\rm img}_s\in{\mathcal S}^{\rm img}_{\bar a,2}$ are guaranteed to be opportunistically stealthy with a probability of at least $50\%$.  

We now have obtained the well-defined ellipsoidal bounded sets ${\mathcal S}^{\rm img}_{\bar a,1}$ and ${\mathcal S}^{\rm img}_{\bar a,2}$, whose structure is the same as that developed for the case of full row-rank $W_a$. Therefore, the corresponding detector design procedure follows a similar spirit to Theorems \ref{prop1}, \ref{theo_reformu} and Algorithm \ref{alg1}, which is not detailed here since no new ideas are required. 

\section{Case Study}\label{case}
In this section, we demonstrate the applicability of the proposed distributionally robust FDI attack detector through a benchmark three-tank system. The laboratory setup known as DTS200 consists of three interconnected tank and two pumps. The detailed parameters of DTS200 can be found in \cite{ding2014data,li2024optimal}. To obtain a discrete-time LTI system, the nonlinear dynamics of DTS200 are linearized around the operating point $x_1=10{\rm cm}$ and $x_3=8{\rm cm}$ with sample time $\Delta T=5$s, resulting in the system \eqref{eq_SS} with state-space matrices
\begin{align*}
   & A = \begin{bmatrix} 
0.8586 & 0.0107 & 0.1304 \\
0.0107 & 0.7958 & 0.1254 \\
0.1303 & 0.1254 & 0.7390 
\end{bmatrix},\\
&B = \begin{bmatrix} 
0.0301 & 0.0001 \\
0.0001 & 0.0290 \\
0.0023 & 0.0022 
\end{bmatrix}, ~
C = I_3,~B_d=I_3,~ D_d=0 .
\end{align*}
The sensor and actuator attack channels are given by 
\begin{align*}
    B_a=\begin{bmatrix}
        I_2&0_{2\times 3}
    \end{bmatrix},  D_a=\begin{bmatrix}
        0_{3\times 2}&I_{3}
    \end{bmatrix}.
\end{align*}
The system is controlled by a feedback controller designed based on an observer and a state-feedback control law. The controller is composed of a Luenberger observer and a state-feedback control law, resulting in the parameters of \eqref{eq_SS_ctrl}: 
\begin{align*}
& A_c = \begin{bmatrix}
0.8078 & 0.0348 & -0.0041 \\
0.0854 & 0.8017 & -0.0307 \\
-0.0583 & 0.0307 & 0.6471
\end{bmatrix},\\
&B_c = \begin{bmatrix}
-1.5368 & -1.7682 & 0.8636 \\
0.8590 & 0.9389 & 0.2613 \\
-1.9897 & -0.6146 & -1.5481
\end{bmatrix},\\
&C_c= \begin{bmatrix}
0.1173 & 0.1297 & -0.2509 \\
-0.3283 & -0.0058 & 0.0179 
\end{bmatrix},~D_c = 0.
\end{align*}

\begin{figure}[htbp]
    \centering
    \includegraphics[width=0.85\linewidth]{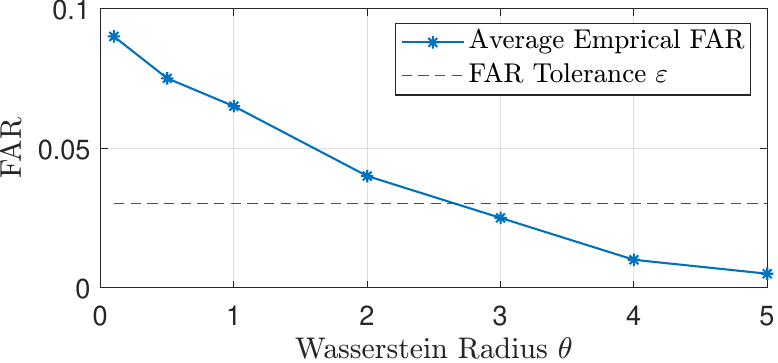}
    \caption{Cross-validation results under different Wasserstein radii $\theta$.}
    \label{fig_radius}
\end{figure}

The residual generator \eqref{eq_fd_parity} with full row-rank $W_a$ and $W_d$ is obtained by using the parity space method \cite{chow1984analytical} with order $s = 4$ is used to generate residuals. Additive disturbances $d(k)$ are generated from a Laplace distribution with covariance $\Sigma_d = 0.01I$ and mean $\mu_d=0$. As for details to implement the proposed method, the convergence rate is chosen as $ \alpha=0.7$ and $N=200$ samples $\{\hat\xi_i\}_{i=1}^{N}$ are used to construct the ambiguity set $\mathcal{D}_{\rm W}(\theta;N)$. The Wasserstein radius $\theta=0.003$ is calibrated via the $5$-fold cross-validation procedure outlined in Remark \ref{rmk_hyper}, as illustrated in Fig. \ref{fig_radius}.
The induced SDPs \eqref{eq_J0}-\eqref{eq_J3} are modeled using YALMIP interface \cite{lofberg2004yalmip} and solved using MOSEK \cite{mosek} in MATLAB R2024a.

\begin{figure}
    \centering
    \includegraphics[width=0.85\linewidth]{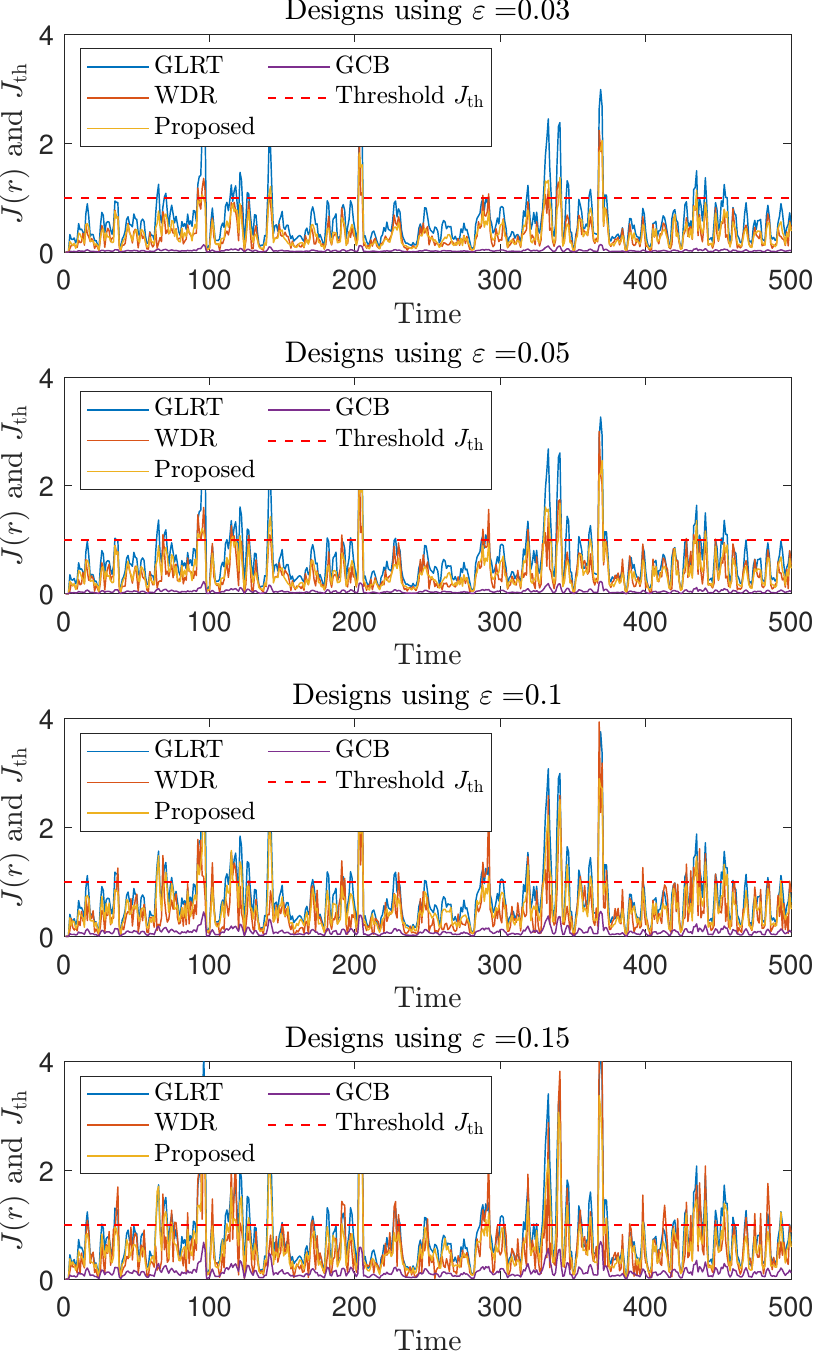}
  \caption{Residual evaluation function $J(r)$ of the proposed method with $\beta=0.7$ and the benchmark designs under the attack-free scenario.}\label{fig_res_free}
\end{figure}

\begin{figure}
    \centering
    \includegraphics[width=0.85\linewidth]{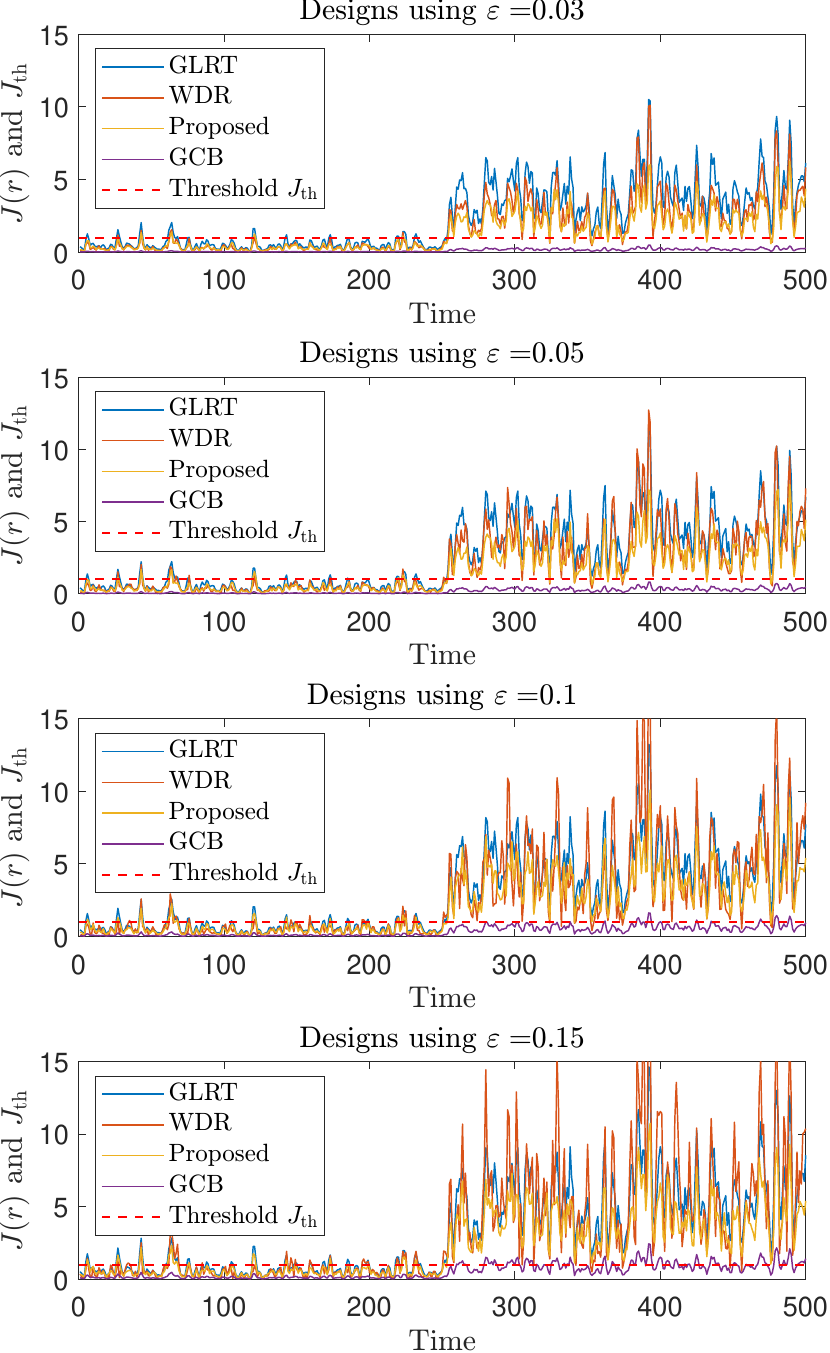}
  \caption{Residual evaluation function $J(r)$ of the proposed method with $\beta=0.7$ and the benchmark designs, where a sequence of uniformly distributed random attacks $\lVert a(k)\rVert_{\infty}\leq0.35$ is injected after the $250$th time step.}\label{fig_res_attack}
\end{figure}

To compare our proposed method with the GLRT design\cite{willsky2003generalized},  the generalized Chebyshev bound (GCB)-based detector \cite{renganathan2020distributionally}, and the WDR detector design \cite{shang2021distributionally} with the same $\theta=0.003$, we performed Monte Carlo simulations with $200,000$ attack and attack-free points. The FARs and attack detection rates (ADRs) of different designs are shown in Table \ref{tab_FAR}, and a representative piece of $500$ samples is shown in Figs. \ref{fig_res_free} and \ref{fig_res_attack} to showcase the evaluation function $J(r)$ under attack-free and attack scenarios. It is seen that the GLRT detector, which is designed under a strict Gaussian assumption, achieves the highest ADRs, but suffers from an "alarm flood" when facing non-Gaussian disturbances. This excessive sensitivity to distributional deviations renders the detector unreliable in practical usage. Meanwhile, the GCB accounts for an excessively broad family of distributions with the same moment information, leading to an overly conservative detector. Although this design guarantees the security of the FAR, it suffers from extremely low ADRs. In contrast, the distributional robustness of WDR and our method enables to maintain the FARs safely below the tolerance level $\varepsilon$, while our design simultaneously achieves both higher ADRs and lower FARs compared with the WDR design \cite{shang2021distributionally}.


\begin{table}[htpb] 
\centering 
\caption{Performance of the proposed detector with $\beta = 0.7$ and benchmark designs evaluated on Monte Carlo simulations.}
\label{tab_FAR}
\begin{tabular}{llcccc} 
\toprule
Method & Index & $\varepsilon=0.03$ & $\varepsilon=0.05$ & $\varepsilon=0.1$ & $\varepsilon=0.15$ \\
\midrule
\multirow{2}{*}{GLRT \cite{willsky2003generalized}} 
    & FAR (\%) &   10.691&   13.254&   18.213&   22.384   \\
    & ADR  (\%) &  99.681&   99.799&   99.912&   99.952    \\
\midrule

\multirow{2}{*}{WDR \cite{shang2021distributionally}} 
    & FAR (\%) &      2.977&    5.498&    9.065&   13.979   \\
    & ADR (\%)  &  95.726&   97.014&   96.904&   97.888   \\
\midrule

\multirow{2}{*}{Proposed}& FAR (\%) &2.426&    4.395&   9.289& 12.183    \\
    & ADR  (\%) & 95.830&   97.931&   99.411&  99.708 \\
\midrule
\multirow{2}{*}{GCB \cite{renganathan2020distributionally}} 
    & FAR (\%) &      0&    0&    0&   0.014   \\
    & ADR (\%)  &  0&   0.012&   8.904&   40.213   \\
\midrule
\bottomrule
\end{tabular}
\end{table}

\begin{figure}
    \centering
    \includegraphics[width=0.95\linewidth]{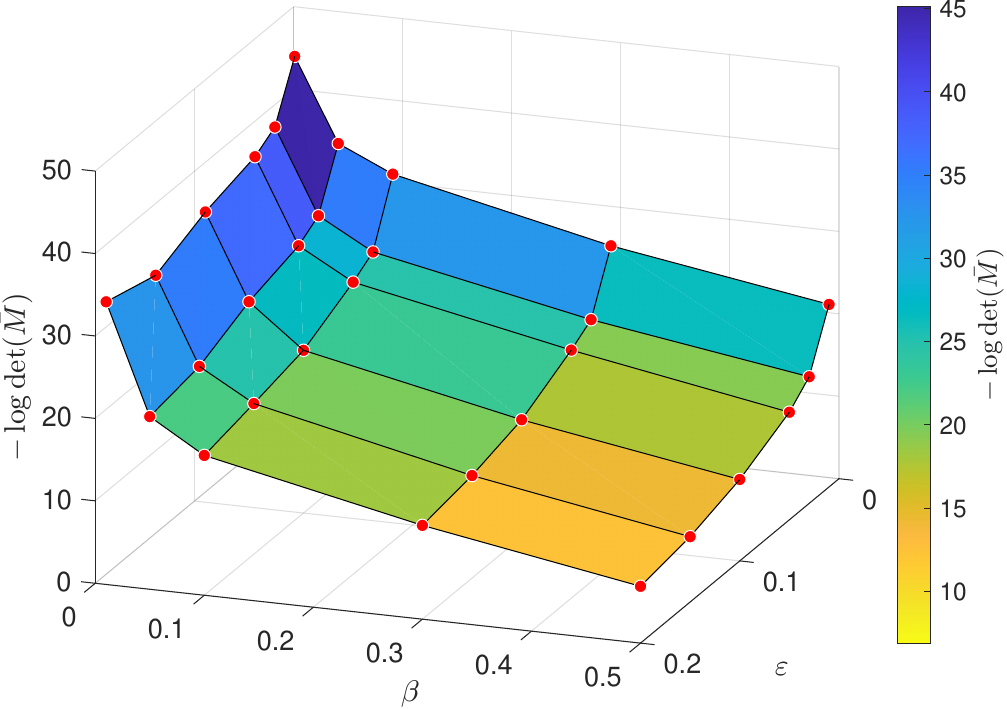}
    \caption{Security metric $\log \det(\bar M)$ of the proposed method versus varying $\varepsilon$ and $\beta$.}\label{fig_logdet}
\end{figure}

\begin{figure}
     \centering
     \includegraphics[width=0.95\linewidth]{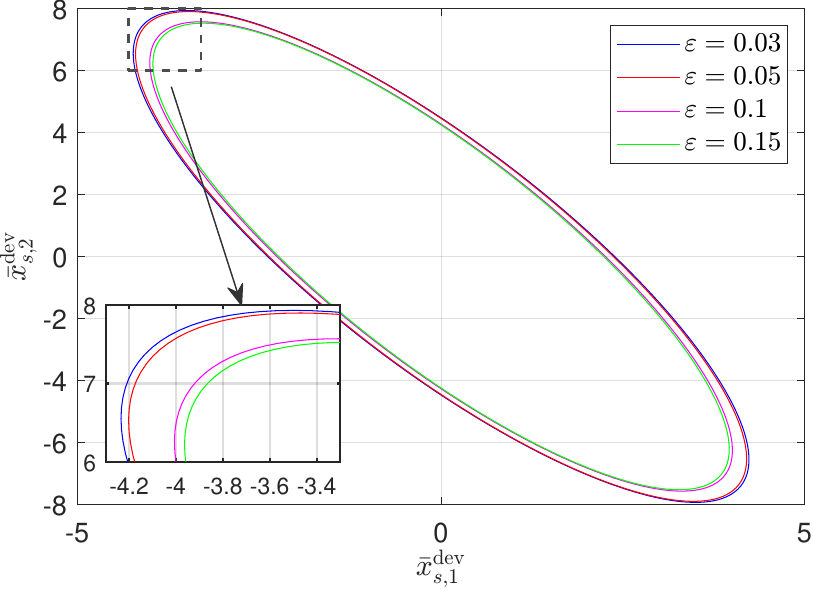}
     \caption{Reachable set approximation of $\bar x^{\rm dev}_{s,[1:2]}$ under $\beta=0.7$ and varying $\varepsilon$.}\label{fig_reach_alpha}
 \end{figure}

\begin{figure}
     \centering
     \includegraphics[width=0.95\linewidth]{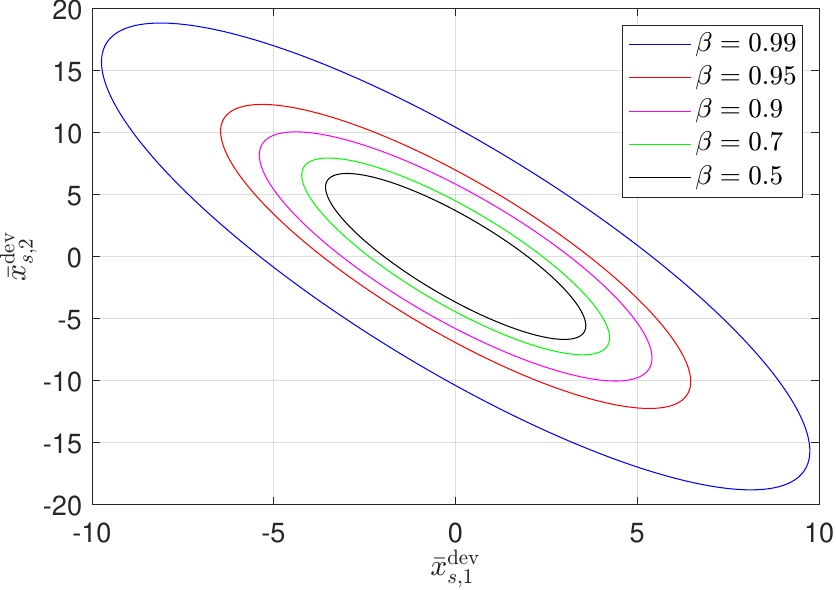}
     \caption{Reachable set approximation of $\bar x^{\rm dev}_{s,[1:2]}$ under $\varepsilon=0.03$ and varying $\beta$.}\label{fig_reach_beta}
 \end{figure}

\begin{figure}
     \centering
     \includegraphics[width=0.95\linewidth]{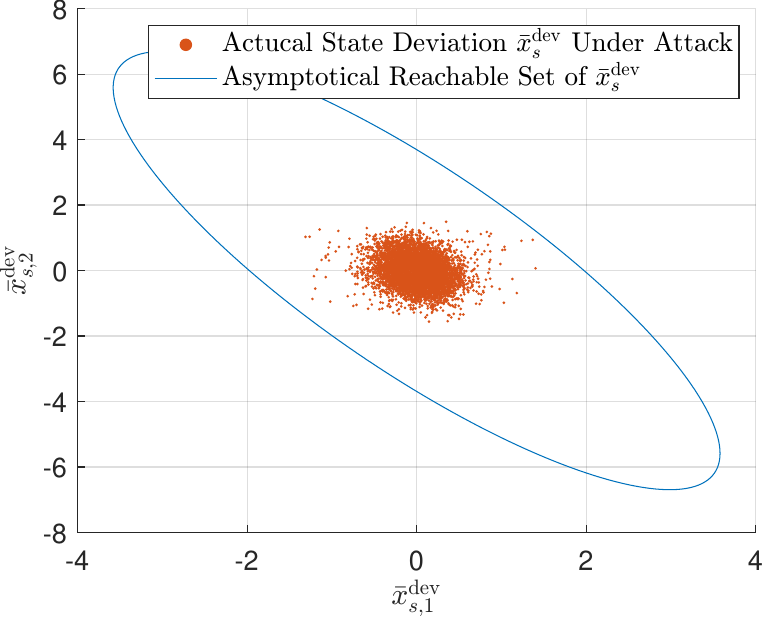}
     \caption{ Predicted asymptotical and actual reachable set of the proposed design with $\beta=0.5$ and $\varepsilon=0.03$.
     }\label{fig_actual}
 \end{figure}

To visualize the high-dimensional ellipsoid $\mathcal{E}^{\infty}_{\bar{x}^{\mathrm{dev}}}$, we project it onto a lower-dimensional subspace, by partitioning $M$ in \eqref{eq_outer_appro} as
\begin{equation*}
    M = \begin{bmatrix}
        M_1 & M_2 \\
        M_2^\top & M_3
    \end{bmatrix},
\end{equation*}
where $M_1 \in \mathbb{R}^{n_{\rm p} \times n_{\rm p}}$. Using the Schur complement, we obtain the following projection  $\mathcal{E}'_{\bar{x}^{\mathrm{dev}}}$ on $\mathbb{R}^{n_{\rm p}}$:
\begin{equation*}
 \mathcal{E}'_{\bar{x}^{\mathrm{dev}}} = \left\{ z \in \mathbb{R}^{n_p} \mid z^\top (M_1 - M_2 M_3^{-1} M_2^\top) z \le 1 \right\}.
\end{equation*}
By setting different values of $\{\varepsilon, \beta\}$, the performance degradation metric $-{ \log\det}(\bar M)$ of the induced detector design is displayed in Fig. \ref{fig_logdet}. The induced projection of the outer approximation ellipsoid  $\mathcal{E}^{\infty}_{\bar{x}^{\mathrm{dev}}}$ onto the first two state dimensions $\bar x^{\rm dev}_{s,1}$ and $\bar x^{\rm dev}_{s,2}$ by using different $\varepsilon$ and $\beta$ is shown in Figs. \ref{fig_reach_alpha} and \ref{fig_reach_beta}. As expected, a general trend is observed that the objective value tends to decrease as $\varepsilon$ increases or $\beta$ decreases. This can be intuitively explained by the balance trade-off in our design. On one hand, enforcing a stricter FAR constraint, i.e., a smaller $\varepsilon$, reduces the sensitivity of the detector, allowing stealthy attacks to cause greater harm and thus resulting in a larger state deviation ellipsoid $\mathcal{E}^{\infty}_{\bar{x}^{\mathrm{dev}}}$. On the other hand, demanding a performance guarantee against a higher confidence set of disturbances, i.e., a larger $\beta$, implies that the performance guarantee must hold for a broader set of disturbances, forcing the design to be more conservative, also reflected in a larger $\mathcal{E}^{\infty}_{\bar{x}^{\mathrm{dev}}}$. To further evaluate the performance, the actual reachable set approximated via $10,000$ Monte Carlo simulations with random initial  $\bar x^{\rm dev}_{s}(0)$ and our predicted asymptotical one under stealthy attack are presented in Fig.\ref{fig_actual}. Clearly, the predicted sets successfully enclose all the state deviation trajectories, thereby validating the effectiveness guarantees of the proposed method.


\begin{figure}
     \centering
     \includegraphics[width=0.95\linewidth]{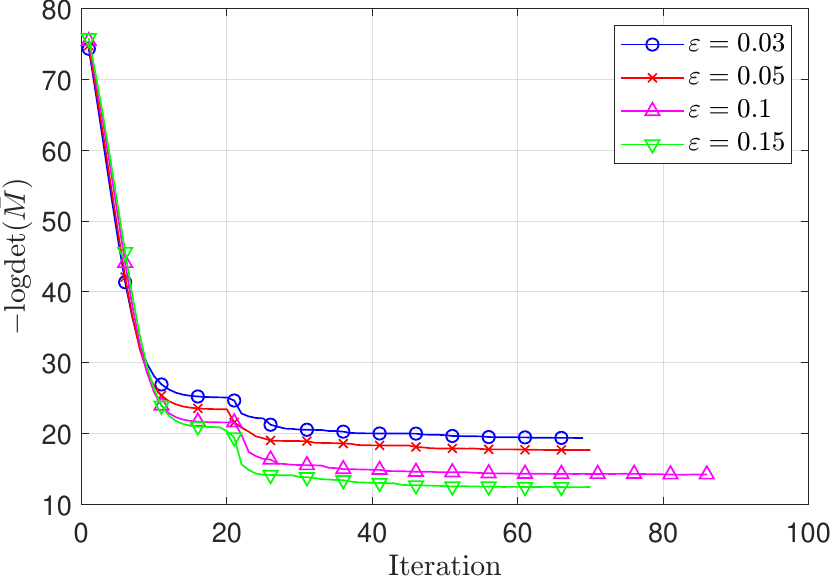}
     \caption{Iterations of the performance degradation metric in Algorithm \ref{alg1} with $\beta = 0.7$ and varying $\varepsilon$.}\label{fig_conv_FAR}
 \end{figure}

\begin{figure}
     \centering
     \includegraphics[width=0.95\linewidth]{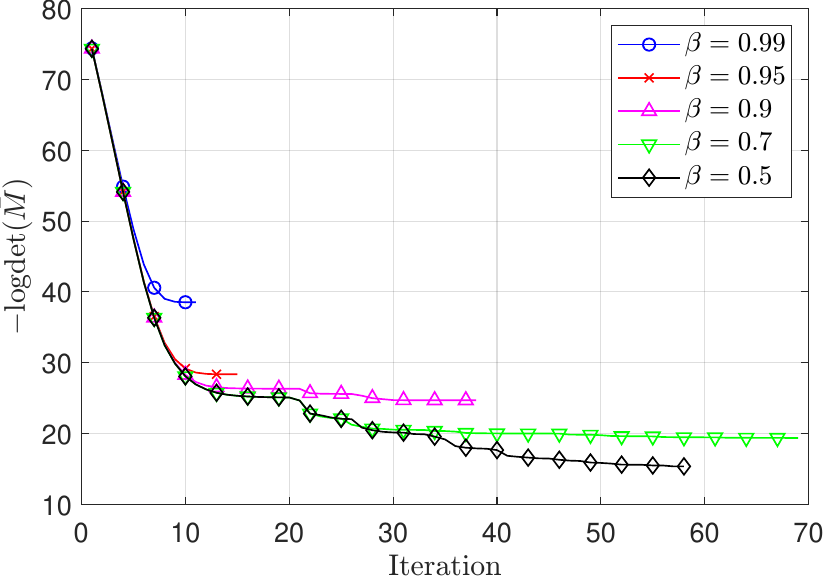}
     \caption{Iterations of the performance degradation metric in Algorithm \ref{alg1} with $\varepsilon = 0.03$ varying $\beta$.}\label{fig_conv_beta}
 \end{figure}

\begin{figure}
    \centering
    \includegraphics[width=0.95\linewidth]{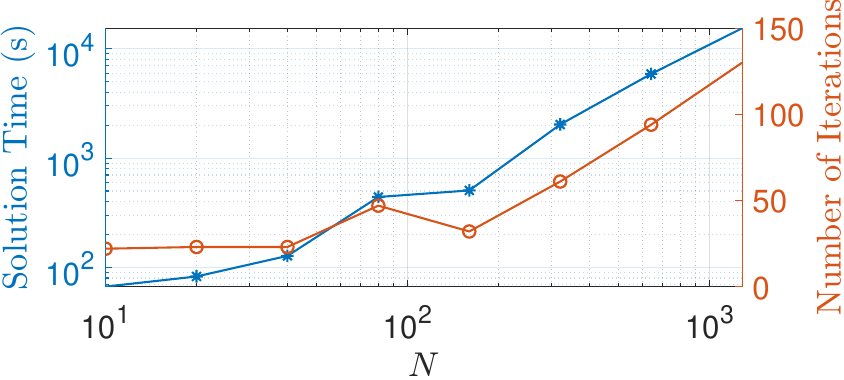}
  \caption{Solution time of Algorithm \ref{alg1} for varying sample sizes $N$.}\label{fig_time}
\end{figure}

Finally, we depict in Figs. \ref{fig_conv_FAR} and \ref{fig_conv_beta} the convergence process of Algorithm 1 under different $\varepsilon$ and $\beta$, which highlights the effectiveness of Algorithm 1. It can be observed that the convergence speed is significantly influenced by $\beta$, whereas it shows insensitivity to $\varepsilon$. Specifically, a larger $\beta$ leads to faster convergence. This is because using a larger $\beta$ reduces the size of the feasible region due to the constraint \eqref{eq_schur}. By operating within a more constrained solution space, the sequential optimization procedure can find a satisfactory solution faster. In contrast, a smaller $\beta$ enlarges the feasible region, leading to a more extensive search and thus a slower convergence. Finally, Fig. \ref{fig_time} illustrates the computational complexity of Algorithm \ref{alg1} with respect to the sample size $N$. It can be observed that the computation time scales significantly with $N$, whereas the number of iterations exhibits a relatively slow increase. Although the proposed design procedure is entirely performed offline, solving SDPs for realistic CPS with hundreds of states remains computationally challenging due to the rapid complexity scaling. Henceforth, using gradient-based optimization to handle these large-scale SDPs is a promising direction for future work.

\section{Conclusion} \label{concl}
In this work, we proposed a novel detector against FDI attacks in the CPS, with a primary focus on the performance degradation caused by stealthy attacks and stochastic disturbance following an unknown distribution. We first presented a distributionally robust performance degradation metric, which is defined by the volume of an asymptotic outer ellipsoidal approximation of the reachable set of state deviation under stealthy attacks. By optimizing this metric, we formulated the detector design problem while using DRCC to control the FAR under a tolerance level, thereby balancing between FAR and security against stealthy attacks.   
Next, to address the intractability of the original design problem involving DRCCs, we reformulate it into a finite-dimensional program with BMIs, and devise a tailored solution algorithm based on sequential minimization to efficiently solve the non-convex problem. Furthermore, we discussed how to extend our detector design scheme to the more general and challenging case of row rank-deficient attack projection matrix, where we formulate the stealthy attack sets in both conservative and optimistic viewpoints. 
Finally, a case study on a three-tank system demonstrated the efficiency of our design and solution algorithm. For future work, a promising direction is to address strictly stealthy attacks under this distributionally robust paradigm.

\bibliographystyle{IEEEtran}
\bibliography{ref}

\end{document}